\documentclass[11pt,twoside,a4paper]{article}
\usepackage{amssymb}
\usepackage{amsmath}
\usepackage[dvips]{graphicx}
\usepackage[all,2cell,ps]{xy}
\title{Embeddability and Stresses of Graphs}
\author{Eran Nevo \footnote{Institute of Mathematics, Hebrew
University of Jerusalem, Jerusalem Israel, E-mail address:
eranevo@math.huji.ac.il. Supported by an I.S.F. grant.}}
 \newtheorem{thm}{Theorem}[section]

 \newtheorem{prop}[thm]{Proposition}
\newtheorem{prob}[thm]{Problem}
\newtheorem{conj}[thm]{Conjecture}

\begin{document}
\maketitle
\begin{abstract}\label{Abstract}
Gluck \cite{Gluck} has proven that triangulated $2$-spheres are
generically $3$-rigid. Equivalently, planar graphs are generically
$3$-stress free.
We show that linklessly embeddable graphs are generically
$4$-stress free. Both of these results are corollaries of the
following theorem: every $K_{r+2}$-minor free graph is generically
$r$-stress free for $1\leq r \leq 4$. (This assertion is false for
$r\geq 6$.) We give an equivalent formulation of this theorem in
the language of symmetric algebraic shifting and show that its
analogue for exterior algebraic shifting also holds. Some further
extensions are detailed.
\end{abstract}
\section{Introduction}\label{Introduction}
Gluck \cite{Gluck} has proven that triangulated $2$-spheres are
generically $3$-rigid. His proof is based on two classical
theorems, of Cauchy and of Steinitz. Cauchy's rigidity theorem
\cite{Cauchy} asserts that any bijection between the vertices of
two (convex) $3$-polytopes which induces a combinatorial
isomorphism, and which induces an isometry of the facets, induces
an isometry of the two polytopes. Gluck actually used Alexandrov's
\cite{Alexandrov} extension of this theorem which relaxes the
condition by replacing the boundaries of the $3$-polytopes with
arbitrary triangulations of them. Steinitz's theorem
\cite{Steinitz} asserts that any polyhedral $2$-sphere is
combinatorially isomorphic to the boundary complex of some
$3$-polytope. It is easy to see that a graph with $n$ vertices and
$3n-6$ edges is generically $3$-rigid iff it is generically
$3$-stress free. Thus, Gluck's theorem can be stated as:
\begin{thm}(Gluck)\label{Gluck}
Planar graphs are generically $3$-stress free.
\end{thm}
We show that also the following relation between embeddability and
rigidity holds:
\begin{thm}\label{Link&Stress}
Linklessly embeddable graphs are generically $4$-stress free.
\end{thm}
Both of these theorems are corollaries of our main theorem:
\begin{thm}\label{Minor&Stress}
For $2\leq r \leq 6$, every $K_{r}$-minor free graph is
generically $(r-2)$-stress free.
\end{thm}
The proof is by induction on the number of vertices, based on
contracting edges possessing a certain property. We make an
essential use of Mader's theorem \cite{Mader} which gives an upper
bound  $(r-2)n-$ $r-1 \choose 2$ on the number of edges in a
$K_{r}$-minor free graph with $n$ vertices, for $r\leq 7$. Indeed,
Theorem \ref{Minor&Stress} can be regarded as a strengthening of
Mader's theorem, as being generically $l$-stress free implies
having at most $ln-$ $l+1 \choose 2$ edges, a fact which is clear
from the equivalent formulation of Theorem \ref{Minor&Stress} in
terms of symmetric algebraic shifting, detailed below. This also
shows that Theorem \ref{Minor&Stress} fails for $r\geq 8$, as is
demonstrated for $r=8$ by $K_{2,2,2,2,2}$, and for $r>8$ by
repeatedly coning over the resulted graph for a smaller $r$ (e.g.
\cite{Song}). It would be interesting to find a proof of Theorem
\ref{Minor&Stress} that avoids using Mader's theorem (and derive
Mader's theorem as a corollary).

Let $\Delta$ denote the algebraic shifting operator, for both
symmetric and exterior versions. The symmetric
case of the following result is equivalent to Theorem \ref{Minor&Stress}:
\begin{thm}\label{mainThm}
The following holds for symmetric and exterior shifting: for every
$2\leq r \leq 6$ and every graph $G$, if $\{r-1,r\} \in \Delta(G)$
then $G$ has a $K_r$ minor.
\end{thm}
As $\Delta(G)$ is shifted (i.e. if $\{a,b\}\in \Delta(G)$ and
$a'\leq a, b'\leq b$ then $\{a',b'\}\in \Delta(G)$) it is
$k$-colorable iff $\{k,k+1\} \notin \Delta(G)$; in this case a
$k$-coloring $f$ would be $f(i)=min\{i,k\}$. Hence, the following
formulation \'{a} la Hadwiger
$$K_{r}\nprec G \Rightarrow \chi(\Delta(G))\leq r-1$$
holds for $r\leq 6$ and is false for $r\geq 8$; the case $r=7$ is
still open. ($\chi(H)$ is the cromaric number of $H$ and $H\nprec
G$ means that $G$ is $H$-minor free.)

\begin{prob}\label{top.minor}
Does Theorem \ref{mainThm} continue to hold when replacing "$K_r$
minor" with "subdivision of $K_r$"?
\end{prob}
The answer is positive for $r=2,3,4$ as in this case $G$ has a
$K_r$ minor iff $G$ has a subdivision of $K_r$ (\cite{Diestel},
Proposition 1.7.2). Mader proved that every graph on $n$ vertices
with more than $3n-6$ edges contains a subdivision of $K_{5}$
\cite{Mader2}. A positive answer in the case $r=5$ would
strengthen this result.

Let $\mu(G)$ denote the Colin de Verdi\`{e}re's parameter of a
graph $G$.
\begin{conj}\label{Colin-conj}
Let $G$ be a graph and let $k$ be a positive integer. If
$\mu(G)\leq k$ then $\{k+1,k+2\} \notin \Delta(G)$.
\end{conj}
For $k=1,2,3,4$ the conjecture holds true. Colin de Verdi\`{e}re
\cite{Colin1} showed that the family $\{G: \mu(G)\leq k\}$ is
closed under taking minors for every $k$. Note that
$\mu(K_{r})=r-1$. By Theorem \ref{mainThm} the conjecture holds
for $k\leq 4$. Another "evidence" is that clique sums do not
violate the conjecture: Suppose that $G_{1}$ and $G_{2}$ satisfy
the conjecture, $G=G_{1}\cup G_{2}$ and  $G_{1}\cap G_{2}$ is a
clique. Let $max\{\mu(G_{1}),\mu(G_{2})\}=k$. By hypothesis,
$\{k+1,k+2\} \notin \Delta(G_{i})$ for $i=1,2$. By \cite{KUL},
Thm.1.2 $\{k+1,k+2\} \notin \Delta(G)$. Also $\mu(G)\geq k$ and
$\Delta(G)$ is shifted, hence $G$ satisfy the conjecture. (Van der
Holst, Lov\'{a}sz and Schrijver \cite{Holst} investigated the
behavior of Colin de Verdi\`{e}re's parameter under taking clique
sums.) Conjecture \ref{Colin-conj} implies $$\mu(G)\leq k
\Rightarrow e\leq kv-(^{k+1}_{\ \ 2})$$ (where $e$ and $v$ are the
numbers of edges and vertices in $G$, respectively) which is not
known either.

This paper is organized as follows: Section \ref{SecRigidity}
provides relevant background in rigidity theory of graphs, Section
\ref{SecMinors} deals with graph minors, in Section \ref{SecProof}
we prove the results about stress freeness mentioned in the
Introduction, Section \ref{SecShifting} deals with algebraic
shifting - both symmetric and exterior, and concludes with a proof
of Theorem \ref{mainThm} and some extensions concerning embeddability into $2$-manifolds.

\section{Rigidity}\label{SecRigidity}
The presentation here is based mainly on Kalai's \cite{LBT-Kalai}.
Let $G=(V,E)$ be a graph. Let $d(a,b)$ denote Euclidian distance
between points $a$ and $b$ in Euclidian space.
A $d$-embedding $f:V\rightarrow \mathbb{R}^{d}$ is called $rigid$
if there exists an $\varepsilon>0$ such that if $g:V\rightarrow
\mathbb{R}^{d}$ satisfies $d(f(v),g(v))<\varepsilon$ for every
$v\in V$ and $d(g(u),g(w))=d(f(u),f(w))$ for every $\{u,w\}\in E$,
then $d(g(u),g(w))=d(f(u),f(w))$ for every $u,w\in V$.
Loosely speaking, $f$ is rigid if any perturbation of it which
preserves the lengths of the edges actually preserves the
distances between any pair of vertices.
$G$ is called $generically\ d-rigid$ if the set of its rigid
$d$-embeddings is open and dense in the topological vector space
of all of its $d$-embeddings.
Given a $d$-embedding $f:V\rightarrow \mathbb{R}^{d}$, a $stress$
w.r.t. $f$ is a function $w:E\rightarrow \mathbb{R}$ s.t. for
every vertex $v\in V$
$$\sum_{u:\{v,u\}\in E}w(\{v,u\})(f(v)-f(u)) =0.$$
$G$ is called $generically$ $d$-$stress$ $free$ if the set of its
$d$-embeddings which has a unique stress ($w=0$) is open and dense
in the space of all of its $d$-embeddings.

Rigidity and stress freeness can be related as follows: Let
$V=[n]$, and let $Rig(G,f)$ be the $dn\times |E|$ matrix
associated with a $d$-embedding $f$ of $V(G)$ defined as follows:
for its column corresponding to $\{v<u\}\in E$ put the vector
$f(v)-f(u)$ (resp. $f(u)-f(v)$) at the entries of the rows
corresponding to $v$ (resp. $u$) and zero otherwise. $G$ is
generically $d$-stress free if $Ker(Rig(G,f))=0$ for a generic $f$
(i.e. for an open dense set of embeddings). $G$ is generically
$d$-rigid if $Im(Rig(G,f))=Im(Rig(K_{V},f)$ for a generic $f$,
where $K_{V}$ is the complete graph on $V=V(G)$. The dimensions of
the kernel and image of $Rig(G,f)$ are independent of the generic
$f$ we choose; we call $R(G)=Rig(G,f)$ the $rigidity\ matrix$ of
$G$.

$Im(Rig(K_{V},f))$ can be described by the following linear
equations:

$(v_{1},..,v_{d})\in \bigoplus_{i=1}^{d}\mathbb{R}^{n}$ belongs to
$Im(Rig(K_{V},f))$ iff
\begin{equation} \label{ImRig1}
\forall 1\leq i \neq j \leq d\ \   <f_{i},v_{j}>=<f_{j},v_{i}>
\end{equation}
\begin{equation}\label{ImRig2}
\forall 1\leq i\leq d\ \  <e,v_{i}>=0
\end{equation}
where $e$ is the all ones vector and $f_{i}$ is the vector of the
$i$th coordinate of the $f(v)$'s, $v\in V$. From this description
it is clear that $rank(Rig(K_{V},f))=dn-$${d+1}\choose{2}$ (see
Asimov and Roth \cite{Asi-Roth} for more details).

 We need the following theorem of Whiteley:
\begin{thm}\label{Whiteley}(Whiteley \cite{Wh})
Let $G'$ be obtained from a graph $G$ by contracting an edge
$\{u,v\}$.

 (a)If $u,v$ have at least $d-1$ common neighbors and
$G'$ is generically $d$-rigid, then $G$ is generically $d$-rigid.

(b)If $u,v$ have at most $d-1$ common neighbors and $G'$ is
generically $d$-stress free, then $G$ is generically $d$-stress
free.
\end{thm}
In Section \ref{SecShifting} we will prove an analogous statement
in the language of exterior shifting. Theorem \ref{Whiteley} gives
an alternative proof of Gluck's theorem (Whiteley \cite{Wh}):
starting with a triangulated $2$-sphere, repeatedly contract edges
with exactly $2$ common neighbors until a tetrahedron is reached
(it is not difficult to show that this is always possible). By
Theorem \ref{Whiteley}(a) it is enough to show that the
tetrahedron is generically $3$-rigid, as is well known (Asimov and
Roth \cite{Asi-Roth}).

For later use, we need the following result about stress-freeness
of a union of graphs.
\begin{thm}\label{weak-clique-sum,Asimov}(Asimov and Roth \cite{Asi-Roth2})
Let $G_{i}=(V_{i},E_{i})$ be $k$-stress free graphs, $i=1,2$ s.t.
$G_{1}\cap G_{2}$ is $k$-rigid. Then $G_{1}\cup G_{2}$ is
$k$-stress free.
\end{thm}

\section{Minors}\label{SecMinors}
All graphs we consider are simple, i.e. with no loops and no
multiple edges. Let $e=\{v,u\}$ be an edge in a graph $G$. By
$conrtacting$ $e$ we mean identifying the vertices $v$ and $u$ and
deleting the loop and one copy of each double edge created by this
identification, to obtain a new (simple) graph. A graph $H$ is
called a $minor$ of a graph $G$, denoted $H\prec G$, if by
repeated contraction of edges we can obtain $H$ from a subgraph of
$G$. In the sequel we shall make an essential use of the following
Theorem of Mader:
\begin{thm}\label{Minor&Edges}(Mader \cite{Mader})
For $3 \leq r \leq 7$, if a graph $G$ on $n$ vertices has no $K_r$
minor then it has at most $(r-2)n-$ $r-1 \choose 2$ edges.
\end{thm}

\begin{prop}\label{anyEdge5}
For $3\leq r \leq 5$:
 If $G$ has an edge and each edge belongs to at least $r-2$ triangles,
 then $G$ has a $K_r$ minor.
\end{prop}
$Proof$: For $r=3$ $G$ actually contains $K_3$ as a subgraph. Let
$G$ have $n$ vetrices and $e$ edges. Assume (by contradiction)
that $K_r\nprec G$. W.l.o.g. $G$ is connected.

For $r=4$, by Theorem \ref{Minor&Edges} $e\leq 2n-3$ hence there
is a vertex $u \in G$ with degree $d(u)\leq 3$. Denote by $N(u)$
the induced subgraph on the neighbors of $u$. For every $v\in
N(u)$, the edge $uv$ belongs to at least two triangles, hence
$N(u)$ is a triangle, and together with $u$ we obtain a $K_4$ as a
subgraph of $G$, a contradiction.

For $r=5$, by Theorem \ref{Minor&Edges} $e\leq 3n-6$ hence there
is a vertex $u \in G$ with degree $d(u)\leq 5$. Also $d(u)\geq 4$
(as we may assume that $u$ is not an isolated vertex). If $d(u)=4$
then the induced subgraph on $\{u\}\cup N(u)$ is $K_{5}$, a
contradiction. Otherwise, $d(u)=5$. Every $v\in N(u)$ has degree
at least $3$ in $N(u)$, hence $e(N(u))\geq \lceil 3\cdot
5/2\rceil=8$. But $K_4\nprec N(u)$, hence $e(N(u))\leq 2\cdot
5-3=7$, a contradiction.$\blacksquare$

\begin{prop}\label{anyEdge6}
If $G$ has an edge and each edge belongs to at least $4$
triangles, then either $G$ has a $K_6$ minor, or $G$ is a clique
sum over $K_r$ for some $r\leq 4$ (i.e. $G=G_{1}\cup G_{2},
G_{1}\cap G_{2}=K_r$, $G_{i}\neq K_r$, $i=1,2$).
\end{prop}
$Proof$: We proceed as in the proof of Proposition \ref{anyEdge5}:
Assume that $K_6\nprec G$. W.l.o.g. $G$ is connected.
By Theorem \ref{Minor&Edges} $e\leq
4n-10$ hence there is a vertex $u \in G$ with degree $d(u)\leq 7$,
also $d(u)\geq 5$. If $d(u)=5$ then $N(u)=K_{5}$, a contradiction.
Actually, since $K_5\nprec N(u)$ and $N(u)$ has at most $7$
vertices each of them of degree at least $4$, Wagner's structure
theorem for $K_5$-minor free graphs (\cite{Diestel}, Theorem
8.3.4) asserts that $N(u)$ is planar.

If $d(u)=6$, then $12=3\cdot 6-6 \geq e(N(u))\geq 4\cdot 6/2=12$
hence $N(u)$ is a triangulation of the $2$-sphere $S^{2}$. If
$d(u)=7$, then $15=3\cdot 7-6 \geq e(N(u))\geq 4\cdot 7/2=14$. We
will show now that $N(u)$ cannot have $14$ edges, hence it is a
triangulation of $S^{2}$: Assume that $N(u)$ has $14$ edges, so
each of its vertices has degree $4$, and $N(u)$ is a triangulation
of $S^{2}$ minus an edge. Let us look on the unique square
(in a planar embedding) and denote its vertices by
$A$. Counting missing edges (there are $7$ of them) shows that
there is one missing edge between the vertices of $N(u)\setminus
A=\{a,b,c\}$, say $\{b,c\}$. we now look at the neighborhood of
$a$ in a planar embedding (it is a $4$-cycle): $b,c$ must be
opposite in this square as $\{b,c\}$ is missing. Hence for $v\in
A\cap N(a)$ we get that $v$ has degree $5$, a contradiction.

Now we are left to deal with the case where $N(u)$ is a
triangulation of $S^{2}$, and hence a maximal $K_5$-minor free
graph. If $G$ is the cone over $N(u)$ with apex $u$, then every
edge in $N(u)$ belongs to at least $3$ triangles in $N(u)$. By
Proposition \ref{anyEdge5}, $N(u)$ has a $K_5$ minor, a
contradiction. Hence there exists a vertex $w\neq u$, $w\in
G\setminus N(u)$.
Denote by $[w]$ the set of all vertices in $G$ connected to $w$ by
a path disjoint from $N(u)$. Denote by $N'(w)$ the induced graph
on the vertices in $N(u)$ that are neighbors of some vertex in
$[w]$. If $N'(w)$ is not a clique, there are two non-neighbors
$x,y\in N'(w)$, and a path through vertices of $[w]$ connecting
them. This path together with the cone over $N(u)$ with apex $u$
form a subgraph of $G$ with a $K_6$ minor, a contradiction.

Suppose $N'(w)$ is a clique (it has at most $4$ vertices, as
$N(u)$ is planar). Then $G$ is a clique sum of two graphs that
strictly contain $N'(w)$: Let $G_{1}$ be the induced graph on
$[w]\cup N'(w)$ and let $G_{2}$ be the induced graph on
$G\setminus [w]$. Then $G=G_{1}\cup G_{2}$ and $G_{1}\cap
G_{2}=N'(w)$. $\blacksquare$

\textbf{Remark}
In view of Theorem \ref{Minor&Edges} for the case $r=7$, we may
expect the following to be true:
\begin{prob}\label{r=7,8}
If $G$ has an edge and each edge belongs to at least $5$
triangles, then either $G$ has a $K_{7}$ minor, or $G$ is a clique
sum over $K_{l}$ for some $l\leq 6$.
\end{prob}
If true, it extends the assertion of Theorem \ref{mainThm} to the
case $r=7$. By now we can show only the weaker assertion
$$\{6,7\}\in\Delta(G) \Rightarrow K_7^-\prec G,$$ using similar
arguments to those used for proving Theorem \ref{mainThm} ($K_7^-$
is $K_7$ minus an edge).
However, if the assertion of Problem \ref{r=7,8} holds for some
$r$, it implies that Theorem \ref{mainThm} holds for this $r$,
hence $e(G)=e(\Delta (G))\leq (r-2)n-$ $r-1 \choose 2$. But as
mentioned in the Introduction, this is false for $r\geq 8$.

\section{Proof of Theorems \ref{Gluck}, \ref{Link&Stress} and \ref{Minor&Stress}}\label{SecProof}
$Proof$ $of$ $Theorem$ $\ref{Minor&Stress}$: For $r=2$ the
assertion of the theorem is trivial. Suppose $K_r\nprec G$, and
contract edges belonging to at most $r-3$ triangles as long as it
is possible. Denote the resulted graph by $G'$. Repeated
application of Theorem \ref{Whiteley} asserts that if $G'$ is
generically $(r-2)$-stress free, then so is $G$. In case $G'$ has
no edges, it is trivially $(r-2)$-stress free. Otherwise, $G'$ has
an edge, and each edge belongs to at least $r-2$ triangles. For
$2<r<6$, by Proposition \ref{anyEdge5} $G'$ has a $K_r$ minor,
hence so has $G$, a contradiction. For $r=6$, by Proposition
\ref{anyEdge6} $G'$ either has a $K_6$ minor which leads to a
contradiction, or $G'$ is a clique sum over $K_r$ for some $r\leq
4$. In the later case, denote $G'=G_{1}\cup G_{2}$, $G_{1}\cap
G_{2}=K_r$. As the graph of a simplex is $k$-rigid for any $k$, by
Theorem \ref{weak-clique-sum,Asimov} it is enough to show that
each $G_{i}$ is generically $(r-2)$-stress free, which follows
from induction hypothesis on the number of vertices. $\blacksquare$

\textbf{Remark} Note that we proved the case $r=5$ without using
Wagner's structure theorem for $K_5$-minor free graphs
(\cite{Diestel}, Theorem 8.3.4), but we used Theorem
\ref{Minor&Edges} of Mader. Alternatively, we can prove the case
$r=5$ avoiding Mader's theorem but using Wagner's theorem and the
'gluing lemma' Theorem \ref{weak-clique-sum,Asimov}. Using
Wagner's structure theorem for $K_{3,3}$-minor free graphs
(\cite{Diestel}, ex.18 on p.185) and Theorem
\ref{weak-clique-sum,Asimov}, we conclude that $K_{3,3}$-minor
free graphs are generically $4$-stress free.

Theorems \ref{Gluck} and \ref{Link&Stress} now follow as easy
corollaries:

$Proof$ $of$ $Theorem$ $\ref{Link&Stress}$: By (the easy part of)
the theorem by Robertson, Seymour and Thomas characterizing
linklessly embeddable graphs by a family of forbidden minors
\cite{Seymour}, a linklessly embeddable graph has no $K_{6}$
minor, hence by Theorem \ref{Minor&Stress} it is generically
$4$-stress free. $\blacksquare$

$Proof$ $of$ $Theorem$ $\ref{Gluck}$: By (the easy part of)
Kuratowski's criterion for planarity of graphs \cite{Kur}, a
planar graph has no $K_{5}$ minor, hence by Theorem
\ref{Minor&Stress} it is generically $3$-stress free. $\blacksquare$

\section{Algebraic shifting}\label{SecShifting}
\subsection{definition of algebraic shifting}
Algebraic shifting is an operator which associates with each
simplicial complex another simplicial complex which is
combinatorially simpler. It was introduced by Kalai \cite{55}. We
follow the definitions and notation of \cite{skira}: Let $K$ be a
simplicial complex on a vertex set $[n]$. The i-th skeleton of $K$
is $K_{i}=\{S\in K: |S|=i+1\}$. For each $1\leq k \leq n$ let
$<_{L}$ be the lexicographic order on $(^{[n]}_{\ k})$, i.e.
$S<_{L}T \Leftrightarrow min\{a:a\in S\triangle T\}\in S$, and let
$\triangleleft_{P}$ be the partial order defined by: Let
$S=\{s_{1}<\dots<s_{k}\}, T=\{t_{1}<\dots<t_{k}\}$,
$S\triangleleft_{P} T$ iff $s_{i}\leq t_{i}$ for every $1\leq
i\leq k$ ($min$ and $\leq$ are taken with respect to the usual
order on $\mathbb{N}$). $K$ is called $shifted$ if
$S\triangleleft_{P} T\in K$ implies $S\in K$.

We now describe exterior shifting: Let $V$ be an $n$-dimensional
vector space over a field $k$ of characteristic zero, with basis
$\{e_{1},\dots,e_{n}\}$. Let $\bigwedge V$ be the graded exterior
algebra over $V$. Denote $e_{S}=e_{s_{1}}\wedge\dots\wedge
e_{s_{j}}$ where $S= \{s_{1}<\dots<s_{j}\}$.
Define the exterior algebra of $K$ by the ring quotient
$$\bigwedge (K)=\bigwedge V/(e_{S}:S\notin K)=\bigwedge V/sp\{e_{S}:S\notin K\}.$$ Let
$\{f_{1},\dots,f_{n}\}$ be a basis of $V$, generic over
$\mathbb{Q}$ with respect to $\{e_{1},\dots,e_{n}\}$, which means
that the entries of the corresponding transition matrix $A$ are
algebraically independent over $\mathbb{Q}$. Let $\tilde{f}_{S}$
be the image of $f_{S}\in \bigwedge V$ in $\bigwedge(K)$. We
choose  a basis for $\bigwedge(K)$ from these images in the greedy
way, to construct the following collection of sets:
$$\Delta^{e}(K)=\bigcup_{i}\{S: \tilde{f}_{S}\notin
sp\{\tilde{f}_{S'}:S'<_{L}S\},|S|=i\}.$$
 The construction is canonic (i.e. independent both of the numbering
 of the vertices of $K$ and of the choice of the generic matrix
 $A$), and results in a shifted simplicial complex.

For symmetric shifting, let us look on the face ring
(Stanley-Reisner ring) of $K$ $k[K]=k[x_{1},..,x_{n}]/I_{K}$ where
$I_{K}$ is the homogenous ideal generated by the monomials whose
support is not in $K$ (grading is by degree). Let
$y_{1},\dots,y_{n}$ be generic linear combinations of
$x_{1},\dots,x_{n}$. We choose a basis for each graded component
of $k[K]$, up to degree $dim(K)+1$, from the canonic projection of the monomials in the
$y_{i}$'s, in the greedy way:
$$GIN(K)=\{m: \tilde{m}\notin sp\{\tilde{m'}:deg(m')=deg(m), m'<_{L}m\}\}$$
(where $\prod y_{i}^{a_{i}}<_{L}\prod y_{i}^{b_{i}}$ iff for
$j=min\{i: a_{i}\neq b_{i}\}$ $a_{j}>b_{j}$). The combinatorial
information in $GIN(K)$ is redundant: if $m\in GIN(K)$ of degree
$i\leq dim(K)$ then $y_{1}m,..,y_{i}m$ are also in $GIN(K)$. Thus,
$GIN(K)$ can be reconstructed from its monomials of the form
$m=y_{i_{1}}\cdot y_{i_{2}}\cdot..\cdot y_{i_{r}}$ where $r\leq
i_{1}\leq i_{2}\leq..\leq i_{r}$, $r\leq dim(K)+1$. Denote this
set by $gin(K)$, and define
$S(m)=\{i_{1}-r+1,i_{2}-r+2,..,i_{r}\}$ for such $m$. The
collection of sets
$$\Delta^{s}(K)=\cup \{S(m): m\in gin(K)\}$$
carries the same combinatorial information as $GIN(K)$. It is a
shifted simplicial complex. Again, the construction is canonic, in
the same sense as for exterior shifting.
\subsection{connection with rigidity and proof of Theorem \ref{mainThm}}\label{sh-rigid}
Let $G$ be a graph. By the results of Lee \cite{Lee},
$\{d+1,d+2\}\notin \Delta^{s}(K)$ iff $G$ is generically
$d$-stress free, as both of these assertions are equivalent to a
zero kernel of the rigidity matrix. We will describe now a similar
statement for exterior shifting in more details; the exterior
analogue of rigidity being Kalai's notion of hyperconnectivity
\cite{56}.

We keep the notation from the previous subsection and follow the
presentation in \cite{56}. Fix $k=\mathbb{R}$. Let $(\bigwedge
V)^{*}\cong\bigwedge(V^{*})$ be the dual of $\bigwedge V$. Fixing
the basis $e=\{e_{1},\dots,e_{n}\}$ induces an inner product on
the degree $j$ part of $\bigwedge V$, denoted $\wedge^{j} V$, for
every $j$: $<f,g>=f^{*}(g)$ is a bilinear extension of
$e_{S}^{*}(e_{T})=\delta_{S,T}$, where $|S|=|T|=j$. Define a left
interior product of $g$ on $f$, where $g,f \in \wedge V$, denoted
$g\lfloor f$, by the requirement:
$$<h,g\lfloor f>=<h\wedge g,f>\ for\ all \ h\in\bigwedge V.$$ Thus,
$g\lfloor f$ is a bilinear function, satisfying
$$e_{T}\lfloor e_{S}=\{^{\pm e_{S\backslash T}\  if\  T\subseteq S} _{0 \ otherwise}$$
where the sign equals $(-1)^{a}$, where $a=|\{(s,t)\in S\times T:
s\notin T, t<s\}|$.

This implies in particular that for $g$ a wedge product of
elements of degree 1, $g \lfloor$ is a boundary operation on
$\bigwedge V$, and in particular on $\bigoplus_i M_{i}(K)$ where
$M_{i}(K)$ is the subspace of $\bigwedge V$ spanned by $\{e_{S}:
S\in K_{i}\}$. Consider the map
$$f(d,i,K): M_{i}(K)\rightarrow \bigoplus_{1}^{d}M_{i-1}(K) \ \ x\mapsto (f_{1}\lfloor
x,...,f_{d}\lfloor x).$$ The dimension of its kernel equals
$|\{S\in \Delta^{e}K: |S|=i+1, S\cap [d]=\emptyset\}|$ (more
details in \cite{KUL}). Kalai \cite{56} called a graph $G$
$d$-$hyperconnected$ if $Im(f(d,1,G))=Im(f(d,1,K_{Ver(G)}))$, and
$d$-$acyclic$ if $Ker(f(d,1,G))=0$. With this terminology, $G$ is
$d$-acyclic iff $\{d+1,d+2\}\notin \Delta^{e}(K)$.

We shall prove now an exterior analogue of Theorem \ref{Whiteley}:
\begin{prop}\label{extWHiteley}
If $G'$ is obtained from $G$ by contracting an edge which belongs
to at most $d-1$ triangles, and $G'$ is $d$-acyclic, then so is
$G$.
\end{prop}
$Proof$: Let $\{v,u\}$ be the edge we contract. Consider the
$dn\times |E|$ matrix $A$ of the map $f(d,1,G)$ w.r.t. the
standard basis, where $f_{i}=\sum_{j=1}^{n}\alpha_{ij}e_{j}$,
$n=|V|$: for its column corresponding to $\{v<u\}\in E$ put the
vector $(\alpha_{1u},..,\alpha_{du})^{T}$ (resp. $-(\alpha_{1v},..,\alpha_{dv})^{T}$)
at the entries of the rows corresponding to $v$ (resp. $u$) and
zero otherwise.

Now replace in $A$ each $\alpha_{iv}$ with $\alpha_{iu}$ to obtain
a new matrix $\hat{A}$. It is enough to show that the columns of
$\hat{A}$ are independent: As the set of $dn\times |E|$ matrices
with independent columns is open (in the Euclidian topology), by
perturbing the $\alpha_{iu}$'s in the places where $\hat{A}$
differs from $A$, we may obtain new generic $\alpha_{iv}$'s
forming a matrix with independent columns. But for every generic
choice of $f_{i}$'s, the map $f(d,1,G)$ has the same rank, hence
we would conclude that the columns of $A$ are independent as well.

Suppose a linear combination of the columns of $\hat{A}$ equals
zero. Let $\bar{A}$ be obtained from $\hat{A}$ by adding the rows
of $v$ to the corresponding rows of $u$, and deleting the rows of
$v$. Thus, the a linear combination with the same coefficients of
the columns of $\bar{A}$ also equals zero. $\bar{A}$ is obtained
from the matrix of $f(d,1,G')$ by adding a zero column (for the
edge $\{v,u\}$) and doubling the columns which corresponds to
common neighbors of $v$ and $u$ in $G$. As $Ker(f(d,1,G'))=0$,
apart from the above mentioned columns the rest have coefficient
zero, and pairs of columns we doubled have opposite sign. Let us
look at the submatrix of $\hat{A}$ consisting of the 'doubled'
columns with vertex $v$ and the column of $\{v,u\}$, restricted to
the rows of $v$: it has generic coefficients, $d$ rows and at most
$d$ columns, hence its columns are independent. Thus, all
coefficients in the above linear combination are zero. $\blacksquare$

We need the following exterior analogue of Theorem
\ref{weak-clique-sum,Asimov}:
\begin{thm}\label{weak-clique-sum,Kalai}(Kalai \cite{56})
Let $G_{i}=(V_{i},E_{i})$ be $k$-acyclic graphs, $i=1,2$ s.t.
$G_{1}\cap G_{2}$ is $k$-hyperconnected. Then $G_{1}\cup G_{2}$ is
$k$-acyclic.
\end{thm}

$Proof\  of\  Theorem\ \ref{mainThm}$: As explained in subsection
\ref{sh-rigid}, Theorem \ref{Minor&Stress} is equivalent to the
symmetric case of Theorem \ref{mainThm}. In the exterior case, the
case $r=2$ is trivial as shifting preserves the $f$-vector. Now we
repeat the proof of Theorem \ref{Minor&Stress} almost word by
word, introducing the following modifications. Replace "Theorem
\ref{Whiteley}" by "Proposition \ref{extWHiteley}". Replace
"stress free" by "acyclic", and "rigid" by "hyperconnected"
everywhere. Replace "Theorem \ref{weak-clique-sum,Asimov}" by
"Theorem \ref{weak-clique-sum,Kalai}". As $G$ is $(r-2)$-acyclic
iff $\{r-1,r\} \notin \Delta^e(G)$, the proof is completed.
$\blacksquare$

\subsection{embedding into $2$-manifolds}\label{surface-sec}
Theorem \ref{Gluck} may be extended to other $2$-manifolds as
follows:
\begin{thm}\label{surface-thm}
Let $M\neq S^2$ be a compact connected $2$-manifold without
boundary, and let $G$ be a graph. Suppose that $\{r-1,r\}\in
\Delta(G)$ and $K_{r}$ can not be embedded in $M$. Then $G$ can
not be embedded in $M$.
\end{thm}
$Proof$: Let $g=g(M)>0$ be the genus of $M$ (e.g. the torus has
genus 1, the projective plane has genus 1/2). Assume by
contradiction that $G$ embeds in $M$. By looking at the rigidity
matrix we note that deleting from $G$ a vertex of degree at most
$r-2$ preserves the existence of $\{r-1,r\}$ in the shifted graph.
Deletion preserves embeddability in $M$ as well. Thus we may
assume that $G$ has minimal degree $\delta(G)\geq r-1$. By Euler
formula $e\leq 3v-6+6g$ (where $e$ and $v$ are the numbers of
edges and vertices in $G$ respectively). Also $e\geq (r-1)v/2$,
hence $v\leq \frac{12g-12}{(r-1)-6}$. Thus
$(r-1)^2-5(r-1)+(6-12g)\leq 0$ which implies $r \leq
(7+\sqrt{1+48g})/2$. But $K_{r}$ can not be embedded in $M$, hence
by Ringel and Youngs \cite{Ringel} proof of Heawood's map-coloring
conjecture $r > (7+\sqrt{1+48g})/2$, a contradiction.$\blacksquare$

\textbf{Remark} For any compact connected $2$-manifold without
boundary of positive genus, $M$, embedded in $\mathbb{R}^3$, two
linked simple closed curves on it exist. One may ask whether the
graph of any triangulated such $M$ is always not linkless. For the
projective plane this is true. It follows from the fact that the
two minimal triangulations of the projective plane (w.r.t. edge
contraction), determined by Barnette \cite{Barnette}, have a minor
from the Petersen family, and hence are not linkless, by the
result of Robertson, Seymour and Thomas \cite{Seymour}. Moreover,
the graph of any polyhedral map of the projective plane is not
linkless, as its $7$ minimal polyhedral maps (w.r.t. edge
contraction), determined by
Barnette \cite{Barnette2}, have graphs equal to $6$ of the members
in Petersen family.

Examining the $21$ minimal triangulations of the torus, see
Lavrenchenko \cite{Lavrenchenko}, we note that $20$ of them have a
$K_6$ minor, and hence are not linkless, but the last one is
linkless, see Figure 1 (one checks that it contains no minor from Petersen's family).
Taking connected sums of this
triangulation, we obtain linkless graphs triangulating any
oriented surface of positive genus. By performing stellar
operations we obtain linkless graphs with arbitrarily many
vertices triangulating any oriented surface of positive genus.

\begin{prob}\label{TP}
Is the graph of a triangulated non orientable $2$-manifold always not linkless?
\end{prob}

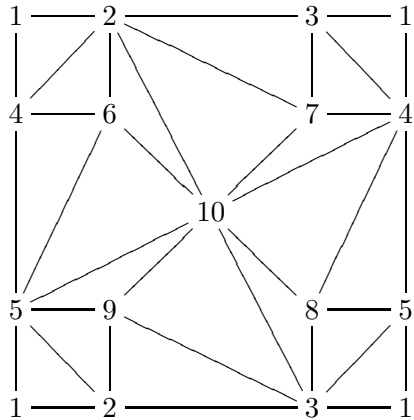
\begin{figure}\label{T21}
\newcommand{\edge}[1]{\ar@{-}[#1]}
\newcommand{\lulab}[1]{\ar@{}[l]^<<{#1}}

\newcommand{\rulab}[1]{\ar@{}[r]^<<{#1}}
\newcommand{\rrulab}[1]{\ar@{}[rr]^<<{#1}}

\newcommand{\ldlab}[1]{\ar@{}[l]^<<{#1}}
\newcommand{\rdlab}[1]{\ar@{}[r]_<<{#1}}
\newcommand{\rrdlab}[1]{\ar@{}[rr]_<<{#1}}
\newcommand{\node}{*+[O][F-]{ }}
\centerline{ \xymatrix{
1  \edge{d} \edge{r} & 2 \edge{d} \edge{rr} \edge{ddr} \edge{rrd} && 3 \edge{d} \edge{r} \edge{dr} & 1 \edge{d}\\
4  \edge{dd} \edge{ur} \edge{r} & 6 \edge{dr} \edge{ddl} && 7 \edge{r} \edge{dl} & 4 \edge{dll} \edge{dd} \edge{ddl} \\
&& 10 \edge{dll} \edge{dl} \edge{ddr} \edge{dr} && \\
5 \edge{d} \edge{r} \edge{dr} & 9 \edge{drr} \edge{d} && 8 \edge{r} \edge{d} & 5 \edge{d} \edge{dl} \\
1 \edge{r} & 2 \edge{rr} && 3 \edge{r} & 1
} } \caption {Linkless graph of a torus}
\end{figure}

\section*{Acknowledgements}
I would like to thank my advisor prof. Gil Kalai for many helpful discussions, and prof. Carsten Thomassen for his contribution to the Minors Section.

\end{document}